\documentclass[aap,MSNbibl,seceqn,dvips]{arximspdf}

\doi{10.1214/12-AAP916} 
\volume{24}
\issue{1}
\pubyear{2014}
\firstpage{114}
\lastpage{130}

\makeatletter
\newcommand{\rrvert}{\vert}
\newcommand{\llvert}{\vert}

\newtheorem{theorem}{Theorem}[section]
\newtheorem{lemma}[theorem]{Lemma}

\makeatother

\begin{document}
\begin{frontmatter}

\title{A Gibbs sampler on the $\lowercase{n}$-simplex}
\runtitle{A Gibbs sampler on the $n$-simplex}

\begin{aug}
\author[A]{\fnms{Aaron} \snm{Smith}\corref{}\thanksref{t1}\ead[label=e1]{asmith3@math.stanford.edu}}
\runauthor{A. Smith}
\affiliation{ICERM, Brown University}
\address[A]{ICERM\\
Brown University\\
121 South Main Street \\
Providence, Rhode Island 02912 \\
USA\\
\printead{e1}} 
\end{aug}

\thankstext{t1}{Supported by a Stanford Graduate Fellowship courtesy of
the Hewlett Foundation.}

\received{\smonth{8} \syear{2011}}
\revised{\smonth{12} \syear{2012}}

%
\begin{abstract}
We determine the mixing time of a simple Gibbs sampler on the unit
simplex, confirming a conjecture of Aldous. The upper bound is based on
a two-step coupling, where the first step is a simple contraction
argument and the second step is a non-Markovian coupling. We also
present a MCMC-based perfect sampling algorithm based on our proof
which can be applied with Gibbs samplers that are harder to analyze.
\end{abstract}

%

\end{frontmatter}

\section{Introduction}\label{sec1}
Given a measure $\mu$ on a convex body $K \subset\mathbb{R}^{n}$, how
can we efficiently obtain independent samples from the distribution of
$\mu$? This problem arises in the computational sciences, and a
frequently-used tool is Markov chain Monte Carlo (MCMC) \cite{Diac08}.
Because MCMC methods produce nearly-independent samples only after a
lengthy mixing period, a long-standing mathematical question is to
analyze the mixing times of the MCMC algorithms in common use.

The analysis of discrete MCMC algorithms is very advanced, with precise
bounds for many difficult problems as well as some general theory that
has received recent exposition in \cite{LPW09,AlFi94}. For samplers on
continuous state spaces, there has been some general theory based on
geometric or coupling arguments (see \cite{LoVe03,Lova98,Yuen01} and
\cite{Rose95}), but many of the techniques built for discrete chains
seem to run into technical difficulties. There are also very few
well-understood simple chains, in stark contrast to the discrete
theory, which has been built on many detailed analyses of specific
chains; though, see \cite{RaWi05a,RaWi05b} for some very nice
analyses of two slower walks on the simplex; \cite{Rose93,Rose94} for
group walks; and \cite{HoJo01} for some applications.
This paper is an attempt to carefully analyze a simple continuous
chain, namely a Gibbs sampler on the $n$-simplex. In addition, it
illustrates the use of two powerful techniques from the discrete
theory: non-Markovian coupling \cite{HaVi03,Borm11,BuKo11} and coupling
from the past \cite{PrWi96}.

The ideas in this paper can be used for a number of other problems. The
analysis was initially motivated by a simpler version of Kac's random
walk on $S(n)$ or $\operatorname{SO}(n)$ (see \cite{Oliv07,Smit12,Jian11} and
especially \cite{Jian12} for recent progress). It is also
a stepping stone toward analysis of Gibbs samplers on more complicated
convex sets, such as contingency tables. In the author's thesis and a
forthcoming note, we use the technique in this paper to improve
existing analyses of these samplers and some others
\cite{Smit11,Smit12}; there is still substantial room for improvement.

In this paper, we will discuss mixing in terms of the popular total
variation distance. For a Markov chain with transition kernel $K$ on a
measurable space $(\Omega, \Sigma)$ and unique stationary distribution
$\pi$, the total variation distance to stationarity after $t$ steps of
a Markov chain started at $\omega\in\Omega$ is given by
\[
\sup_{A \in\Sigma} \bigl| K^{t}(\omega, A) - \pi(A) \bigr|.
\]

Most of this paper will be concerned with a specific Gibbs sampler
$X_{t}$ on the $n$-simplex $\Delta_{n} = \{ X \in\mathbb{R}^{n} |
\sum_{i=1}^{n} X[i] = 1; X[i] \geq0 \}$ whose stationary distribution is
the uniform distribution on $\Delta_{n}$. To take a move in this Markov
chain, begin by choosing $1 \leq i < j \leq n$ and $\lambda\in[0,1]$
independently and uniformly. Then set
%
%
\begin{eqnarray}
\label{EqMoveRep} X_{t+1}[i] &=& \lambda\bigl(X_{t}[i] +
X_{t}[j]\bigr),
\nonumber
\\
X_{t+1}[j] &=& (1 - \lambda) \bigl(X_{t}[i] +
X_{t}[j]\bigr),
\\
X_{t+1}[k] &=& X_{t}[k] \qquad (k \neq i,j).
\nonumber
\end{eqnarray}
This sampler was first mentioned in \cite{AlFi94}, where the mixing
time was shown to be $O(n^{2} \log n )$. Aldous suggested in his list
of open problems that the correct mixing time was $O(n \log n )$, and
we confirm this, also demonstrating a pre-cutoff window of moderate size:

%
\begin{theorem}[(Simplex mixing time)]\label{theo1.1}
Fix $C> 3$ and $n$ satisfying $n > \max(4096, 2C + \frac{7}{2})$ and
$\frac{n}{\log n} > \frac{3 ({1/2} + C)C}{{C/2} -
{1/4}} $. If $K_{n}^{t}$ is the $t$-step transition kernel for
the Gibbs sampler described above, and $U_{n}$ is the uniform
distribution on $\Delta_{n}$, then for all $t > 10 C n \log n $, $x
\in\Delta_{n}$ and $A \subset\Delta_{n}$ measurable,
\[
\bigl| K_{n}^{t}(x, A) - U_{n}(A) \bigr|<
n^{3-C} + 2 n^{-{C/2} -
{1/4}} + 4 n^{{11/4} - C}.
\]
On the other hand, for $0 < C < 1$ and $t < (1 - C) n \log n$,
\[
\liminf_{n \rightarrow\infty} \sup_{x \in\Delta_{n}} \sup
_{A \subset
\Delta_{n}} \bigl| K_{n}^{t}(x,A) -
U_{n}(A) \bigr|= 1.
\]
\end{theorem}

The conditions on the constant $C$ are not onerous. Choosing $C = 4$
gives a mixing time of at most $40 n \log n $ that is effective for $n
> 4096$.

Sections \ref{sec2}--\ref{sec4} are devoted to proving the upper bound
of Theorem \ref{theo1.1},
and Section \ref{sec5} proves the lower bound. In Section \ref{sec6},
we briefly discuss
applications of our method to closely related Markov chains. In
Section \ref{sec7}, we use the ideas of the proof to develop a perfect sampling
algorithm with wider applicability.

\section{Notation, basic lemmas and strategy}\label{sec2}
We recall that a coupling of Markov chains with transition kernel $K$
is a process $(X_{t}, Y_{t})$ so that marginally, both $X_{t}$ and
$Y_{t}$ are Markov chains with transition kernel $K$. The proof relies
on the following standard lemma (see \cite{LPW09}, Theorem 5.2---they
work in discrete space, but their proof does not rely on this assumption):
%
%
\begin{lemma}[(Fundamental coupling lemma)]\label{lem2.1}
Assume $(X_{t}, Y_{t})$ is a coupling of Markov chains such that if
$X_{s} = Y_{s}$, then $X_{t} = Y_{t}$ for all $t > s$. Assume also that
$X_{0} = x$ and $Y_{0}$ is distributed according to the stationary
distribution of $K$. Define the random time $\tau$ to be the first time
at which $X_{t} = Y_{t}$. Then $\sup_{A \in\Sigma} | K^{t}(x,A) -
\pi(A) |\leq P[\tau> t]$.
\end{lemma}

Throughout this note, we are interested in a coupling of Markov chains
$(X_{t}, Y_{t})$, where $X_{0}$ starts according to some distribution
of our choosing,
$Y_{0}$ starts out uniformly over the simplex and both marginally
evolve as the Gibbs sampler being studied. We will describe a joint
evolution of our two chains $X_{t}$ and $Y_{t}$, such that at a
specific time, the probability of having coupled is very high. The
method for proving this is slightly unusual. In most coupling proofs,
including the non-Markovian coupling in \cite{HaVi03}, there is an
attempt to make the two chains get closer throughout the process. In
our method, we attempt to couple only at a specific final time, and
include many moves that are likely to increase the distance between the
chains by a large amount. In fact, our joint distribution will
generally assign 0 probability to coupling at any prior time.

In order to develop our global joint coupling, we describe two possible
one-step couplings of $X_{t}$ and $Y_{t}$. These are the
``proportional'' coupling and the ``subset'' coupling. Throughout, we
will always choose to update entries at the same coordinates $i,j$ in
both $X_{t}$ and $Y_{t}$ at every step; only the uniform variable
$\lambda$ used in representation (\ref{EqMoveRep}) sometimes differs.
Because of this, we often describe the couplings by describing only how
the update variables $\lambda$ are coupled.

In the proportional coupling, we choose an $i,j$ and $\lambda$ for
$Y_{t}$, and then use the same choices for $X_{t}$ in representation
(\ref{EqMoveRep}), so that, for example, entry $i$ in $Y_{t}$ is
updated to $\lambda(Y_{t}[i] + Y_{t}[j])$ while entry $i$ in $X_{t}$ is
updated to $\lambda(X_{t}[i] + X_{t}[j])$. The subset coupling is
slightly more complicated. As before, we choose two coordinates $i,j$
to be updated in both chains. Next, define the weight $w(S,X)$ that a
vector $X$ gives to a subset $S \subset[n]$ to be $w(S,X) = \sum_{s
\in S} X[s]$. A subset coupling of $X_{t}$ and $Y_{t}$ will always be
with respect to some specific subset $S \subset[n] = \{ 1, 2,\ldots,
n \}$. If either $i,j \in S$ or $i,j \notin S$, perform a proportional
coupling. Otherwise, assume without loss of generality that $i \in S$
and $j \notin S$ and also that $X_{t}[i] + X_{t}[j] \geq Y_{t}[i] +
Y_{t}[j]$. In this case, call a coupling of $X_{t+1}$ and $Y_{t+1}$
conditioned on $X_{t}, Y_{t}, i$ and $j$ a subset coupling if
\[
P\bigl[w(X_{t+1},S) = w(Y_{t+1},S)\bigr] \geq
\frac{Y_{t}[i] + Y_{t}[j] - |
\sum_{k \in S / \{ i \} } ( Y_{t}[k] - X_{t}[k] ) |
}{X_{t}[i] + X_{t}[j]}.
\]

We will say a subset coupling has succeeded if $w(X_{t+1},S) =
w(Y_{t+1},S)$, and that it has failed otherwise. We will generally not
be concerned with what happens when a subset coupling has failed. We
will check now that such a coupling exists. Note that, conditioned on
$X_{t}$ and the coordinates $i,j$ updated at time $t$, the weight
$w(S,X_{t+1})$ is uniformly distributed on $[\sum_{k \in S / \{i \}}
X_{t}[k], \sum_{k \in S / \{i \}} X_{t}[k] + X_{t}[i] + X_{t}[j]]$.
Similarly, conditioned on $Y_{t},i$ and $j$, $w(S,Y_{t+1})$ is
uniformly distributed on $[\sum_{k \in S / \{i \}} X_{t}[k], \sum_{k
\in S / \{i \}} X_{t}[k] + X_{t}[i] + X_{t}[j]]$.

%
\begin{lemma}[(Total variation distance of two uniform distributions)]
\label{LemmaFromReviewer}\label{lem2.2}
Let $U$ be distributed uniformly on $[a, a+b]$ and let $U'$ be
distributed uniformly on $[a', a' + b']$. Assume $b \leq b'$. Then
$\|\mathcal{L}(U) - \mathcal{L}(U') \|_{\mathrm{TV}} \leq
1 -
\frac{b - | a - a' |}{b'}$.
\end{lemma}

\begin{pf}
Note that $U$ has density $f(x) = \frac{1}{b} \mathbf{1}_{x \in[a,
a+b]}$, and $U'$ has density $g(x) = \frac{1}{b'} \mathbf{1}_{x \in
[a', a'+b']}$. Thus, the total variation distance between them is given by
\begin{eqnarray*}
\bigl\|\mathcal{L}(U) - \mathcal{L}\bigl(U'\bigr) \bigr\|_{\mathrm{TV}} &=& 1 - \int_{x} \min\bigl(f(x), g(x)
\bigr) \,dx
\\
&=& 1 - \frac{1}{b'} \int_{x \in[a,a+b] \cap[a', a' + b']} 1 \,dx
\\
&=& 1 - \frac{1}{b'} \bigl[\min\bigl(a+b,a' + b'
\bigr) - \max\bigl(a,a'\bigr)\bigr]
\\
&\leq& 1 - \frac{1}{b'} \bigl[b + \min\bigl(a,a'\bigr) - \max
\bigl(a,a'\bigr)\bigr]
\\
&=& 1 - \frac{b - | a - a' |}{b'}.
\end{eqnarray*}
\upqed\end{pf}

Since it is always possible to couple two random variables $W,Z$ so
that $P[Z = W] = 1 - \|\mathcal{L}(Z) - \mathcal{L}(W)
\|_{\mathrm{TV}}$, Lemma \ref{LemmaFromReviewer} implies that subset
couplings exist.

We now give a rough and nonrigorous description of the proof strategy,
which proceeds by describing a two-step coupling of $X_{t}$ and
$Y_{t}$. For the first $T_{1}$ steps, $X_{t}$ and $Y_{t}$ evolve always
under the proportional coupling. This coupling is Markovian, and we
prove that under this coupling, the two chains are very close in
sup-norm with high probability after about $n \log n$ steps. In the
next phase, we record the updated coordinates $(i(t), j(t))$ from time
$T_{1}$ until a specified time $T = T_{1} + T_{2}$. This information is
used to construct a nested sequence $P_{t}$ of partitions of the set of
coordinates $[n]$. With high probability, the sequence will satisfy
$P_{T_{1}} = \{ [n] \}$ and $P_{T_{1} + T_{2}} = \{ \{ 1\}, \{ 2
\},\ldots, \{ n \} \}$. We will then couple $X_{t}$ to $Y_{t}$ step by
step, using only information about the future that is contained in
$P_{t}$, using a proportional coupling for some steps and a subset
coupling for others. We then show that it is possible to keep
$w(S,X_{t}) = w(S,Y_{t})$ for all $S \in P_{t}$ with high probability.
If all of these high-probability events occur, then the final partition
consists of only singletons, and this implies that $X_{T}[i] =
Y_{T}[i]$ for $1 \leq i \leq n$. The two main difficulties are
constructing the partition and showing that $X_{t}$ and $Y_{t}$ remain
close throughout the second phase.

It is worth pointing out that the dependence of the coupling on the
future is in fact necessary to get the correct mixing time, or indeed
any bound that is $o(n^{2})$. This is analogous to the well-known fact
that no Markovian coupling of the random transposition walk on $S_{n}$
can give a coupling time that is $o(n^{2})$. See Lemma 8 of \cite
{Borm11} for a short proof of this fact for the walk on $S_{n}$ which
applies essentially without modification to this Gibbs sampler.

Here is a list of some commonly used variables that have been reserved,
for reference while reading:\\[9pt]
$X_{t}$, the Markov chain of interest.
\\
$Y_{t}$, another instance of the Markov chain, started at stationarity.
\\
$P_{t}$, a set partition of $[n]$.
\\
$S$, a piece of a partition.
\\
$i,j$, coordinates we update.
\\
$\lambda, \lambda_{x}, \lambda_{y}$, uniform random variable
used to update a chain, or chains $X_{t}$ and $Y_{t}$.
\\
$w(S,X)$, the weight assigned by vector $X$ to a subset $S
\subset[n]$.

\section{First coupling stage}\label{sec3}
Define $Z_{t} = \| X_{t} - Y_{t} \|_{2}^{2}$. The
following provides an upper bound for $E[Z_{t}]$ under the proportional
coupling described above:

%
\begin{lemma}[(Burn-in)]\label{lem3.1}
Let $X_{t}$ and $Y_{t}$ be two copies of the Markov chain coupled by
the proportional coupling, and $Z_{t}$ defined as above. After $s
\geq\frac{3}{2} d n \log n$ steps of the proportional coupling,
$E[Z_{s}] \leq2n^{-d}$.
\end{lemma}

\begin{pf} The proof is by a one-step contraction estimate. Assume
$X_{t}$ and $Y_{t}$ are coupled by the proportional coupling from time
0 onwards. Let $\mathcal{F}_{t}$ be the $\sigma$-algebra generated by
the random variables $X_{t}$ and $Y_{t}$; note that $Z_{t}$ is
$\mathcal
{F}_{t}$-measurable. Then, under the proportional coupling, the
following equality comes from conditioning on the coordinates $(i,j) =
(i(t),j(t))$ updated at time $t$,
\begin{eqnarray*}
&&
E[Z_{t+1} |\mathcal{F}_{t}] \\
&&\qquad= \frac{1}{n(n-1)} \sum
_{i \neq j} E \biggl[ \lambda^{2}
\bigl(X_{t}[i] + X_{t}[j] - Y_{t}[i] -
Y_{t}[j]\bigr)^{2}
\\
&&\qquad\quad\hspace*{71pt}{}+ (1 - \lambda)^{2} \bigl(X_{t}[i] + X_{t}[j] -
Y_{t}[i] - Y_{t}[j]\bigr)^{2}\\
&&\qquad\quad\hspace*{71pt}\hspace*{89.5pt}{} + \sum
_{k \neq i,j} \bigl(X_{t}[k] - Y_{t}[k]
\bigr)^{2} \biggr].
\end{eqnarray*}
Note $E[\lambda^{2}] = E[(1- \lambda)^{2}] = \frac{1}{3}$. Expanding
the above, we obtain
\begin{eqnarray*}
&&
E[Z_{t+1} |\mathcal{F}_{t}] \\
&&\qquad= \frac{1}{n(n-1)} \sum
_{i \neq j} \biggl[ \frac{2}{3}
\bigl(X_{t}[i] - Y_{t}[i]\bigr)^{2}
+ \frac{2}{3}\bigl(X_{t}[j] - Y_{t}[j]
\bigr)^{2}\\
&&\qquad\quad\hspace*{60.5pt}{} + \frac{4}{3} \bigl(X_{t}[i] -
Y_{t}[i]\bigr) \bigl(X_{t}[j] - Y_{t}[j]\bigr) \\
&&\qquad\quad\hspace*{130pt}{} +
\sum_{k \neq i,j} \bigl(X_{t}[k] -
Y_{t}[k]\bigr)^{2} \biggr].
\end{eqnarray*}
Collecting coefficients of $Z_{t}$ and using the fact that $Z_{t} =
\sum_{k} (X_{t}[k] - Y_{t}[k])^{2}$, this equals
\[
\biggl( 1 - \frac{2}{3n} \biggr) Z_{t} + \frac{4}{3n(n-1)} \sum
_{i \neq
j} \bigl(X_{t}[i] - Y_{t}[i]
\bigr) \bigl(X_{t}[j] - Y_{t}[j]\bigr).
\]

Noting that $\sum_{i=1}^{n} (X_{t}[i] - Y_{t}[i]) = 0$, the last term
can be rewritten as
\begin{eqnarray*}
\sum_{i \neq j} \bigl(X_{t}[i] -
Y_{t}[i]\bigr) \bigl(X_{t}[j] - Y_{t}[j]\bigr) &=&
\Biggl( \sum_{i=1}^{n} \bigl(X_{t}[i]
- Y_{t}[i]\bigr) \Biggr)^{2} - \sum
_{i=1}^{n} \bigl(X_{t}[i] -
Y_{t}[i]\bigr)^{2}
\\
&=& -Z_{t}.
\end{eqnarray*}

Putting this together, we find that
\[
E[Z_{t+1} |\mathcal{F}_{t}] = \biggl( 1 -
\frac{2}{3n} - \frac
{4}{3n(n-1)} \biggr) Z_{t}.
\]
And so in particular,
\begin{eqnarray*}
E[Z_{t} |\mathcal{F}_{0}] &=& E\bigl[E[Z_{t}
|\mathcal{F}_{t-1}] |\mathcal{F}_{0}\bigr]
\\
&\leq&\biggl( 1 - \frac{2}{3n} \biggr) E[Z_{t-1} |
\mathcal{F}_{0}].
\end{eqnarray*}
By induction on $t$, it is then easy to see that
\[
E[Z_{t} |\mathcal{F}_{0}] \leq\biggl(1 -
\frac{2}{3n} \biggr)^{t} Z_{0}.
\]

Bound $Z_{0}$ by
\begin{eqnarray*}
Z_{0} &\leq& \sum_{k} \bigl(
X_{0}[k]^{2} + Y_{0}[k]^{2} \bigr)
\\
&\leq& \sum_{k} X_{0}[k] +
Y_{0}[k]
\\
&=& 2.
\end{eqnarray*}
We conclude that at times $t \geq\frac{3}{2} d n \log n$, $E[Z_{t}]
\leq2 n^{-d}$.
\end{pf}
Using the\vspace*{1pt} obvious inequality $| X_{t}[i] - Y_{t}[i] |\leq
\sqrt{Z_{t}}$ and Markov's inequality, $P[| X_{t}[i] - Y_{t}[i]
|>
\delta] \leq\delta^{-1} n^{-d/2}$ for all $\delta> 0$ and
$i \in[n]$.

\section{Second coupling stage}\label{sec4}

Let $T = (\frac{1}{2} + \varepsilon) n \log n$ be fixed, for some
$\varepsilon> 0$ to be decided later. Let $Y_{0}$ be chosen from the
uniform distribution on the simplex, and let $X_{0}$ satisfy $\| X_{0}
- Y_{0} \|_{1} \leq n^{-d}$. We describe a coupling $(X_{t}, Y_{t})$
from time $0$ to time $T$ with the property that $X_{T} = Y_{T}$ with
high probability as $n$ goes to infinity, for any fixed $\varepsilon>
0$ and $d$ sufficiently large. First, we choose a sequence of pairs of
distinct elements $1 \leq i(t) \neq j(t) \leq n$ independently and
uniformly for times $0 \leq t \leq T$. These pairs $(i(t), j(t))$ will
be the coordinates updated at time $t$ in both $X_{t}$ and $Y_{t}$.
Then define a sequence of graphs $G_{t}$ for $0 \leq t \leq T-1$ to
have vertex set $[n]$ and edge set $E_{t} = \{ (i(t), j(t)), (i(t+1),
j(t+1)),\ldots, (i(T-1), j(T-1)) \}$, throwing out repeated edges, if
any. We also define $G_{T}$ to be the graph on $[n]$ with no edges.
From this sequence, construct a sequence of partitions of $[n]$, $P(0),
P(1),\ldots, P(T)$ by letting the sets in $P_{t}$ be exactly the
connected components of $G_{t}$.

Since the edges satisfy $E_{s} \subset E_{t}$ for every $s > t$, it is
clear that for any $A \in P_{s}$, there must be some $B \in P_{t}$ with
$A \subset B$. In this sense, the sequence of partitions is nested.
Also note from the construction that either $P_{t}$ and $P_{t+1}$ are
the same, or they differ by having a single set in $P_{t}$ split into
two sets in $P_{t+1}$. Define the sequence of marked time $0 \leq t_{1}
< \cdots< t_{k} = T -1$ as the times at which $P_{t_{\ell}} \neq
P_{t_{\ell}+1}$. Then, for marked time $t_{\ell}$, define $S(t_{\ell
},1)$ and $S(t_{\ell},2)$ to be the two sets that were split apart at
time $t_{\ell}$, labeled so that $| S(t_{\ell}, 1) |\leq
|
S(t_{\ell}, 2) |$. Note that there are at most $n-1$ marked times,
and that there are $n-1$ if and only if $P_{0} = [n]$.

Note also that $P_{0} = [n]$ if and only if $G_{0}$ is connected. The
question of whether or not the random graph $G_{0}$ is connected is a
classical question in random graph theory. The following result, found
in \cite{Boll01} among other places, is good enough for our purposes:

%
\begin{lemma}[(Connectedness for Erdos--Renyi graphs)]\label{lem4.1}
Let $\varepsilon> 0$ be fixed, and let $T = T_{\varepsilon}$ be the first
time that $( \frac{1}{2} + \varepsilon) n \log n$ distinct edges have
been chosen. Then the probability that $G_{0}$ is connected is at least
$1 - n^{-\varepsilon}$.
\end{lemma}

This has the immediate corollary:

%
\begin{lemma}[(Connectedness for $G_{0}$)]\label{lem4.2}
Let $\varepsilon> 0$, and assume $n> 4$ satisfies $\frac{n}{\log n} >
\frac{3 (1 + 2 \varepsilon)({1/2} + 2 \varepsilon)}{\varepsilon
}$. Then
let $T > ( \frac{1}{2} + 2 \varepsilon) n \log n$. Then the probability
that $G_{0}$ is connected is at least $1 - 2n^{-\varepsilon}$.
\end{lemma}

\begin{pf}
Ignoring the ordering of vertices in edges, define
\[
A_{t} = \mathbf{1}_{(i(t),j(t)) \notin\{(i(0),j(0)),\ldots,
(i(t-1),j(t-1)) \}} \mathbf{1}_{t-1 < T_{\varepsilon}}.
\]
Note that $P[A_{t} = 1 | A_{1},\ldots, A_{t-1}] \leq\frac{(1 + 2
\varepsilon) \log n }{n-1}$. In particular, if $B$ has binomial $((\frac
{1}{2} + 2 \varepsilon) n \log n, \frac{2({1/2} + \varepsilon)
\log n
}{n-1})$ distribution, then $P[\sum_{s = 0}^{({1/2} + 2
\varepsilon)
n \log n} A_{s} > x] \leq P[B > x]$ for all $x>0$. For $n$ satisfying
$\frac{n}{\log n} > \frac{3 \varepsilon}{(1 + 2 \varepsilon)({1/2}
+ 2
\varepsilon)}$, Chernoff's inequality gives the bound
\begin{eqnarray*}
P \bigl[T_{\varepsilon} > \bigl( \tfrac{1}{2} + 2 \varepsilon\bigr) n
\log n \bigr] &\leq& P[B > \varepsilon n \log n ]
\\
&\leq& e^{-{n \varepsilon}/{2}},
\end{eqnarray*}
which is less than $n^{-\varepsilon}$ for $n \geq4$. Let $\mathcal
{E}_{T}$ be the event that $G_{0}$ is disconnected. Since $P[G_{0}]
\leq P[T_{\varepsilon} > T] + P[E_{T} | T_{\varepsilon} < T]$, the result
follows immediately from this bound on $T_{\varepsilon}$ and Lemma \ref{lem4.2}.
\end{pf}

Having constructed this partition, we now couple $X_{t}$ and $Y_{t}$
for time \mbox{$0 \leq t \leq T$}. First, we need to choose the coordinates to
update; we do this by updating coordinates $i(t)$ and $j(t)$ at time
$t$ in both chains. Next, we must describe the coupling of the
coordinates. If $t$ is a marked time, then perform a subset coupling
for the set $S(t,1)$. Otherwise, do a proportional coupling. Note that
if $t$ is a marked time, then one of $i(t)$ or $j(t)$ is in $S(t,1)$
and the other is in $S(t,2)$, so the coupling proceeds according to the
description in Section \ref{sec2} in the case $i \in S$, $j \in S^{c}$.
We claim
that this couples the two walks by time $T$ with high probability:

%
\begin{lemma}[(Coupling for close chains)]\label{lem4.3}
For $\varepsilon> 0$, $d > \frac{11}{2}$, $n > \max(d + \frac{3}{2},
4096)$, $\frac{n}{\log n } > \frac{3 (1 + 2 \varepsilon)({1/2} + 2
\varepsilon)}{\varepsilon}$ and $T > ( \frac{1}{2} + 2 \varepsilon)
n \log n $, the coupling described in this section has the property
\[
P[X_{T} \neq Y_{T}] \leq2n^{-\varepsilon} +
5n^{({15-2d})/{4}}.
\]
\end{lemma}

We begin by showing that subset couplings succeed with high probability:

%
\begin{lemma}[(Subset coupling)]\label{lem4.4}
Assume $n \geq6$, and let $(X_{t},Y_{t})$ be a pair of elements of
$\Delta_{n}$ satisfying $\sup_{k} | X_{t}[k] - Y_{t}[k] |=
n^{-f}$ and $\inf_{k} X_{t}[k],\break  \inf_{k} Y_{t}[k] \geq2 n^{-b}$, with
$f \geq b + 1$. Then for all $S \subset[n]$ and update coordinates $i
\in S$, $j \notin S$, $P[w(X_{t+1},S) = w(Y_{t+1},S)] \geq1 - 3n^{b + 1
- f}$ under the subset coupling.
\end{lemma}
\begin{pf}
Assume that $X_{t}[i] + X_{t}[j] \geq Y_{t}[i] + Y_{t}[j]$. Then, from
its definition, the subset coupling succeeds with probability at least
\begin{eqnarray*}
&&
P\bigl[w(X_{t+1},S) = w(Y_{t+1},S)\bigr] \\
&&\qquad\geq
\frac{Y_{t}[i] + Y_{t}[j] - |
\sum_{k \in S / \{ i \} } ( Y_{t}[k] - X_{t}[k] ) |}{X_{t}[i] +
X_{t}[j]}
\\
&&\qquad\geq \frac{Y_{t}[i] + Y_{t}[j] - 2 | S | n^{-f}}{Y_{t}[i] +
Y_{t}[j] + 4 n^{-f}}
\\
&&\qquad\geq \bigl(1 - 2n^{1-f+b}\bigr) \bigl(1 - 4 n^{-f + b} - 8
n^{-2f + 2b}\bigr),
\end{eqnarray*}
which, for $n \geq6$, is at least $1 - 3n^{b + 1 - f}$.
\end{pf}

Having bounded the probability of failure when $X_{t}$, $Y_{t}$ are
close, we must show that they remain close as long as all subset
couplings succeed. For $S \subset[n]$, define $\| X \|_{S} = \sum_{s \in S} | X[s] |$. Then:

%
\begin{lemma}[(Closeness)]\label{lem4.5}
Let $X_{t}, Y_{t}$ be coupled as described above, and assume that
$P_{0} = \{ [n] \}$, that all subset couplings up to time $t$ have
succeeded and that $\| X_{0} - Y_{0} \|_{1} <
\varepsilon$. Then $\| X_{t} - Y_{t} \|_{S} <
\varepsilon$
for every $S$ in $P_{t}$
\end{lemma}
\begin{pf}
There are two types of coupling to take care of. For a proportional
coupling with coordinates $i$ and $j$,
\begin{eqnarray*}
&& \bigl\llvert X_{t+1}[i] - Y_{t+1}[i] \bigr\rrvert+ \bigl
\llvert X_{t+1}[j] - Y_{t+1}[j] \bigr\rrvert
\\
&&\qquad= \lambda_{t} \bigl\llvert X_{t}[i] + X_{t}[j] -
Y_{t}[i] - Y_{t}[j] \bigr\rrvert\\
&&\qquad\quad{}+ (1 - \lambda_{t})
\bigl\llvert X_{t}[i] + X_{t}[j] - Y_{t}[i] -
Y_{t}[j] \bigr\rrvert
\\
&&\qquad\leq\bigl\llvert X_{t}[i] - Y_{t}[i] \bigr\rrvert+ \bigl
\llvert X_{t}[j] - Y_{t}[j] \bigr\rrvert.
\end{eqnarray*}
Since $i$ and $j$ always connect elements of the same set in $P_{t}$,
this shows that proportional couplings never increase $\|
X_{t} - Y_{t} \|_{S}$.
Otherwise, assume that at time $t$ we had a successful subset coupling
for the subset $S$ along edges $i$ and $j$. Without loss of generality,
assume that $i \in S:= S(t,1)$ and $j \in R:
= S(t,2)$. Since $w(X_{0}, [n]) = w(Y_{0}, [n]) = 1$, and all subset
couplings up to time $t$ have succeeded, we have $w(X_{t}, Q) =
w(Y_{t}, Q)$ for all $Q \in P_{t}$. In particular, $w(X_{t}, S \cup R)
= w(Y_{t}, S \cup R)$. Then we note that
\begin{eqnarray*}
X_{t+1}[i] - Y_{t+1}[i] &=& \sum_{s \in S \setminus\{ i \} }
\bigl(Y_{t}[s] - X_{t}[s]\bigr)
\\
&=& X_{t}[i] - Y_{t}[i] + \sum_{s \in R}
\bigl(X_{t}[s] - Y_{t}[s]\bigr)
\end{eqnarray*}
and so
\[
\bigl\llvert X_{t+1}[i] - Y_{t+1}[i] \bigr\rrvert\leq\bigl
\llvert X_{t}[i] - Y_{t}[i] \bigr\rrvert+ \|
X_{t} - Y_{t} \|_{R},
\]
which immediately implies that
\[
\| X_{t+1} - Y_{t+1} \|_{S} \leq
\| X_{t} - Y_{t} \|_{S \cup R}.
\]

An analogous calculation shows that
\[
\| X_{t+1} - Y_{t+1} \|_{R} \leq
\| X_{t} - Y_{t} \|_{R \cup S}
\]
as well. By induction on $t$, this implies that $\| X_{t+1} -
Y_{t+1}\|_{S} \leq\| X_{0} - Y_{0} \|_{1}$
and $\| X_{t+1} - Y_{t+1}\|_{R} \leq\|
X_{0} - Y_{0} \|_{1}$.
\end{pf}

%
\begin{lemma}[(Largeness)] \label{LemmaLargeness}
$P[\inf_{1 \leq i \leq n} \inf_{0 \leq t \leq n^{2}-1} Y_{t}[i] \leq
n^{-4.5 -k}] \leq2 n^{-k}$ for $n > \max(2k, 4096)$.
\end{lemma}
\begin{pf}
Let $q_{1},\ldots, q_{n}$ be independent random variables chosen from
the exponential distribution with mean 1, and let $Q = \sum_{i=1}^{n}
q_{n}$. It is well known (see, e.g., Algorithm 2.7.1 of \cite{RuMe98})
that $ ( \frac{q_{1}}{Q},\ldots, \frac{q_{n}}{Q} )$ is
distributed\vspace*{-2pt} uniformly on the simplex $\Delta_{n}$. In particular,
$Y_{t} \stackrel{D}{=} ( \frac{q_{1}}{Q},\ldots, \frac{q_{n}}{Q}
)$. Taking a union bound over $1 \leq i \leq n$ and $0 \leq t
\leq n^{2}-1$, it is thus sufficient to show
\[
P \biggl[ \frac{q_{1}}{Q} \leq n^{-1.5 - k} \biggr] \leq2
n^{-k}.
\]

Let $E$ be the event that $\frac{q_{1}}{Q} < n^{-1.5 - k}$, $E_{1}$ the
event that $q_{1} < n^{-k - 0.25}$, and $E_{2}$ the event that $Q >
n^{1.25}$, and observe that $E \subset E_{1} \cup E_{2}$. It is
immediate that
\begin{eqnarray*}
P[E_{1}] &=& 1 - e^{-n^{-k - 0.25}}
\\
&\leq& n^{-k - 0.25} + \tfrac{1}{2} n^{-2k - 0.5}.
\end{eqnarray*}

For $n > 4096$, this is certainly less than $n^{-k}$. To bound the
probability that $Q$ is large, note that for all $0 < \theta< 1$,
\begin{eqnarray*}
E\bigl[e^{\theta Q}\bigr] &=& E\bigl[e^{\theta q_{1}}\bigr]^{n}
\\
&=& \frac{1}{(1-\theta)^{n}}.
\end{eqnarray*}

Setting $\theta= 1 - n^{-0.25}$, Markov's inequality gives
\[
P[E_{2}] \leq e^{({1/4}) n \log n + n - n^{1.25}}.
\]

It is straightforward to check that, for $n > \max(2k, 4096)$, this is
less than $n^{-k}$. Since $P[E] \leq P[E_{1}] + P[E_{2}]$, this proves
the lemma.
\end{pf}

Finally, it is possible to prove that in fact most couplings will succeed:

%
\begin{lemma}[(Weight lemma)]\label{lem4.7}
Fix $d > \frac{11}{2}$ and $n > \max(k + \frac{3}{2}, 4096)$. Assume
$P_{0} = \{ [n] \}$ and that $\| X_{0} - Y_{0} \|_{1}
\leq n^{-d}$. Let $E$ be the event that the equality
%
%
\begin{equation}
\label{EqWeightsWork} w(S,X_{t}) = w(S,Y_{t})
\end{equation}
holds for all $0 \leq t \leq T$ and all $S \subset P(t)$. Then
\[
P[E] \geq1 - 5n^{({15 - 2d})/{4}}.
\]
\end{lemma}
\begin{pf} The equality (\ref{EqWeightsWork}) clearly holds at time $0$.
Also note that if it holds at an unmarked time $t$, it must also hold
at time $t+1$, since at unmarked times the weights of parts $S$ of the
partition $P_{t}$ cannot change in either $X_{t}$ or $Y_{t}$. So,
assume that equality (\ref{EqWeightsWork}) holds for all times $t \leq
t_{k}$ for some marked time $t_{k}$. If the subset coupling is
successful at time $t_{k}$, then $w(S(t_{k}+1,1),X_{t_{k}+1}) =
w(S(t_{k}+1,1),Y_{t_{k}+1})$ by construction. However, by the
assumption that equality (\ref{EqWeightsWork}) holds until time $t_{k}$,
$w(S(t_{k},1) \cup S(t_{k},2), X_{t_{k}}) = w(S(t_{k},1) \cup
S(t_{k},2), Y_{t_{k}})$. Since $w(A \cup B, X) = w(A,X) + w(B,X)$ for
any disjoint sets $A,B$ and any vector $X$, this implies
$w(S(t_{k}+12),X_{t_{k}+1}) = w(S(t_{k}+1,2),Y_{t_{k}+1})$ as well.
Since none of the other parts of $P_{t_{k}}$ change weight, this
implies that $w(S, X_{t_{k}+1}) = w(S,Y_{t_{k}+1})$ holds for all $S
\in P_{t_{k}+1}$.

It remains to bound only the probability that the first subset coupling
to fail occurs at time $t_{k}$. By Lemma \ref{lem4.3} and the
assumption of this lemma,
%
%
\begin{eqnarray}
\label{IneqSizeBd} \| X_{t_{k}} - Y_{t_{k}} \|_{1} &\leq& \sum_{S \in
P_{t_{k}}} \|
X_{t_{k}} - Y_{t_{k}} \|_{S}
\nonumber
\\
&\leq& \sum_{S \in P_{t_{k}}} n^{-d}
\\
&\leq& n^{1-d}.\nonumber
\end{eqnarray}

Set $q = \frac{d}{2} + \frac{3}{4}$. By Lemma \ref{LemmaLargeness},
$\inf_{i,t} Y_{t}[i] \geq n^{-q}$ with probability at least $1 - 2
n^{4.5 - q}$. Assuming this holds, Lemma \ref{lem4.4} along with inequality
(\ref{IneqSizeBd}) implies that any particular subset coupling succeeds
with probability at least $1 - 3 n^{2+q-d}$. Taking a union bound over
all at most $n-1$ subset couplings, all subset couplings succeed with
probability at least $1 - 3n^{3+q-d} - 2 n^{4.5-q} = 1 - 5n^{({15 -
2d})/{4}}$.
\end{pf}

It is now time to prove Lemma \ref{lem4.2}. Recall that if $P_{0} =
[n]$ and all
components $Q \in P_{t}$ satisfy $w(Q, X_{t}) = w(Q, Y_{t})$ for $0
\leq t \leq T$, then at time $T$ the two walks have coupled. There are
only two ways for this to fail to happen. The first is the event
$E_{1}$ that $P_{0} \neq[n]$. By Lemma \ref{lem4.2}, $P[E_{1}] \leq2n^{-
\varepsilon}$. The second is the event $E_{2}$ that at least one subset
coupling fails. By\vspace*{1pt} Lemma \ref{lem4.7} and our assumption that $\frac
{n}{\log n }
> \frac{1}{2} + 2 \varepsilon$, which implies $T < n^{2}$, we have the
bound $P[E_{2}] \leq5 n^{({15 - 2d})/{4}}$. Combining these two
bounds proves the lemma.

Finally, we prove Theorem \ref{theo1.1}. We will run the proportional
coupling until time $T_{1} = 9 Cn \log n$, and then we will run the
second phase coupling from time $T_{1}$ until time $T = 10 C n \log n
$. There\vspace*{1pt} are only two ways to have $X_{T} \neq Y_{T}$. The
first is the event $E_{1}$ that $\| X_{T_{1}} - Y_{T_{1}} \|_{1} >
n^{-(2C+2)}$. By the\vspace*{1pt} comment immediately after Lemma
\ref{lem3.1}, $P[E_{1}] \leq n^{3-C}$. The second is the event $E_{2}$
that the second phase coupling fails. By Lemma \ref{lem4.3}, $P[E_{2}]
\leq2 n^{{C}/{2} - {1}/{4}} + 4 n^{{11/4} - C}$. Combining these two
bounds proves the theorem.

We also note that it is possible to improve the top of the pre-cutoff
window from 30 to 12 by being more careful in the above proofs, but
there is no hope of actually proving a cutoff without a substantially
new argument.

\section{Lower bound}\label{sec5}
Since our walk is over a continuous space, the total variation distance
to stationarity of the Markov chain at time $t$ must be at least the
probability that not all coordinates have been chosen by time $t$.
Since only two coordinates are chosen at a time, the classical
coupon-collector results in \cite{ErRe61} tell us that at time $T =
\frac{1}{2} n (\log n - c)$, $\sup_{A \in\Sigma} | K_{n}^{T}(x, A)
- \pi(A) |\geq1 - \exp(-\exp(c)) + o(1)$ as $n$ goes to infinity.

It is possible to improve the constant a little bit. Let $X_{0} =
(1,0,\ldots,0)$, and let $Q_{t} \in\{0, 1 \}^{n}$ be a vector keeping
track of updates in $X_{t}$, started at $Q_{0} = (0,0,\ldots, 0)$. If
coordinates $i$ and $j$ are updated in $X_{t}$ at time $t$, set
$Q_{t+1}[i] = Q_{t+1}[j] = 1$ if at least one of $X_{t}[i]$, $X_{t}[j]$
are nonzero, and set $Q_{t+1}[k] = Q_{t}[k]$ for all $k \neq i,j$. If
$X_{t}[i] = X_{t}[j] = 0$, then set $Q_{t+1} = Q_{t}$. Next, let $\tau
_{j} = \inf\{t | Q_{t+1} \neq Q_{t}, t > \tau_{j-1} \}$ with
$\tau
_{0} = 0$. We note that $E[\tau_{1}] = \frac{n}{2}$, and for $j>1$,
$E[\tau_{j}] = \frac{n(n-1)}{2j(n-j)}$. Thus, letting $\tau= \sum
_{j=1}^{n-1} \tau_{j}$,
\begin{eqnarray*}
E[\tau] &=& \frac{n}{2} + \frac{n^{2}}{2} \sum
_{j=2}^{n-1} \frac
{1}{j(n-j)}
\\
&=& n \log n + O(n).
\end{eqnarray*}
Similarly, since $\tau_{i}$ and $\tau_{j}$ are independent for $i
\neq
j$, it is easy to calculate that the variance $V[\tau] \leq6 n^{2}$.
By Chebyshev's inequality, for all $\varepsilon> 0$ and $n$
sufficiently large,
%
%
\begin{equation}
\label{IneqChebLowerBound} P\bigl[\tau< (1 - \varepsilon) n \log
n\bigr] = O
\biggl( \frac{1}{\log n^{2}} \biggr).
\end{equation}
Finally, observe that for $t < \tau$, at least one entry of $X_{t}$ is
0, and so taking $H_{j} = \{ X \in\Delta_{n} | X[j] = 0 \}$ and $A
\in\Sigma$ to be $\bigcup_{j} H_{j}$, we find $|
K_{n}^{T}((1,0,\ldots,\break 0), A) - U_{n}(A) |\geq P[T < \tau]$.
Combined with inequality (\ref{IneqChebLowerBound}), this proves the
lower bound on the mixing time.

\section{Closely related walks}\label{sec6}
It is worth pointing out a small number of cases where the above
argument goes
through with very few changes. The first allows us to go from sampling
from the
uniform distribution to sampling from a large class of distributions on
the simplex, including symmetric Dirichlet distributions. At each step
of the random walk, instead of choosing $\lambda$ according to the
uniform distribution on $[0, 1]$, choose it according to some other
distribution with twice differentiable cdf $F$ satisfying $F[x] =
1-F[1-x]$ for all $0 \leq x \leq\frac{1}{2}$. Then the above arguments
show that the total mixing time is
$O(n \log n \frac{\| F'' \|_{\infty} + 1}{1 - 2
E[\lambda^{2}]} )$, essentially without modifcation.

It is also possible to apply this argument to the discrete analogue of
the simplex, in which $M$ indistinguishable balls are stored in n
boxes; these are known as $M$-compositions of $n$. The analogous Markov
chain involves choosing two boxes, holding $N$ balls between them, at
every step, and putting $0 \leq k \leq N$ of them in the first box with
probability $\frac{1}{N+1}$, and the remainder in the second box. The
arguments given above apply to the discrete chain, giving a mixing
bound of order $O(n \log n)$, but there need to be enough balls for the
continuous approximation to be good at each step. A
straightforward\vspace*{1pt} step-through of the argument gives a bound
of $O(n \log n)$ for $M > n^{18.5}$ above. Aldous' greedy argument,
which gives an upper bound of $O(n^{2} \log n)$, holds for $M >
n^{5.5}$.

The follow-up paper \cite{Smit11} will discuss a wider variety of
related walks, requiring larger modifcations.

\section{Perfect sampling on the simplex}\label{sec7}
In this section, we discuss how the two-chain coupling described above
can be modified into a grand coupling, and how to use this fact to
create a perfect sampling algorithm. Before describing the algorithm,
we mention that it is not a practical way to obtain uniform points on
the simplex. However, the same algorithm can be used to obtain samples
from the other distributions on the simplex mentioned in Section \ref{sec6},
many of which are a priori much harder to sample from. The method is
also of some interest as a relatively rare instance of a coupling from
the past (CFTP) algorithm which does not use monotonicity or
anti-monotonicity.

To begin, we recall the CFTP algorithm, described in greater detail in
\cite{PrWi96}. First, choose some large time $T$, and start a copy of
the Markov chain $X_{-T}^{\omega}$ for each $\omega$ in the sample
space $\Omega$. Next, couple all of the chains from time $-T$ to time
0. If the chains have coalesced by time 0, the resulting single value
is distributed according to the stationary distribution of the chain.
If not, we couple chains started at all points from $-2T$ to $T$ and
keep the evolution from $-T$ to 0, then from $-3T$ to $-2T$ keeping the
evolution from $-2T$ to 0, and so on until coalescence at 0 has
occurred.

For Markov chains on a finite state space, it is easy in theory to
construct a grand coupling that will eventually coalesce, though bad
couplings are very inefficient. In practice, even on finite chains,
CFTP is only used if the chain has some very special properties. The
most popular properties are monotonicity and its twin antimonotonicity.
Briefly, we introduce a partial order $\leq$ on $\Omega$, and say that
a coupling of two chains $X_{t}$, $Y_{t}$ is monotone if $X_{0} \leq
Y_{0}$ implies $X_{t} \leq Y_{t}$ for all $t > 0$. It is then easy to
see that if our grand coupling is monotone, it is sufficient to keep
track of chains started at maximal and minimal elements of the poset.
If they have coupled, all states have coupled.

For Markov chains on infinite state spaces, many grand couplings will
never coalesce, and of course we cannot keep track of all of the
starting values on a computer. Some chains have a monotonicity
property, but such a property is not obvious for the simplex model.
Despite this, there is a fairly efficient perfect sampling algorithm
that requires tracking only $n+1$ points (and a little extra overhead
each time an epoch of length $T$ fails to coalesce).\looseness=-1

Let $X_{t}^{v}$ be a copy of the Markov chain started at $v = (v[1],
v[2],\ldots, v[n])$ at time $0$, and let $e_{j}$ be the $j$th standard
unit basis vector. We construct a grand coupling of the chains
$X_{t}^{v}$ as follows. For time $0 < t < T_{1}$, do a proportional
coupling. That is, at each time $t$, choose coordinates $i(t), j(t)$
and parameter~$\lambda(t)$, and update all chains using these three
numbers. We claim that for each~$t$, there exists a matrix $M_{t}[i,j]$
such that for any $v$, $X_{t}^{v}[i] = \sum_{j=1}^{n} M_{t}[i,j] v[j]$.
To see this, observe that $X_{t+1} = M_{i(t),j(t),\lambda(t)} X_{t}$,
where $M_{i(t),j(t),\lambda(t)}[i(t),i(t)] = M_{i(t),j(t),\lambda
(t)}[i(t),j(t)] = \lambda(t)$, $M_{i(t),j(t),\lambda(t)}[j(t),i(t)] =
M_{i(t),j(t),\lambda(t)}[j(t), j(t)] = (1- \lambda(t))$,
$M_{i(t),j(t),\lambda(t)}[k,k] = 1$ for $k \notin\{ i(t), j(t) \}$ and
all other entries are 0. We can then write $M_{t}[i,j] = \prod_{s < t}
M_{i(s),j(s),\lambda(s)}$. Thus, for $v, w \in\Delta_{n}$,
\begin{eqnarray*}
\bigl\| X_{t}^{v} - X_{t}^{w} \bigr\|_{1} &=& \sum_{i=1}^{n} \bigl|
X_{t}^{v}[i] - X_{t}^{w}[i] \bigr|
\\
&=& \sum_{i=1}^{n} \Biggl|\sum
_{j=1}^{n} M_{t}[i,j] \bigl(v[j] - w[j]
\bigr) \Biggr|
\\
&\leq& \sum_{i=1}^{n} n \sup
_{j} M_{t}[i,j]
\\
&\leq& \sum_{i=1}^{n} n \sum
_{j,k} \bigl| X_{t}^{e_{j}} [i] -
X_{t}^{e_{k}} [i] \bigr|,
\end{eqnarray*}
which gives
%
%
\begin{equation}
\label{IneqCftp} \bigl\| X_{t}^{v} -
X_{t}^{w} \bigr\|_{1} \leq n \sum
_{j,k} \bigl\| X_{t}^{e_{j}} -
X_{t}^{e_{k}} \bigr\|_{1}.
\end{equation}
Applying Lemma \ref{lem3.1} to the expectation of $\|
X_{t}^{e_{j}} -
X_{t}^{e_{k}} \|_{2}$ for all distinct pairs $j,k$, using
Markov's inequality to bound the probability\vspace*{1pt} that the $L^{2}$ norm is
large, and finally noting that $\| X_{t}^{e_{j}} -
X_{t}^{e_{k}} \|_{1} \leq n \| X_{t}^{e_{j}} -
X_{t}^{e_{k}} \|_{2}$, we find that for $t > \frac{3}{2} d n
\log n$,
\[
P\bigl[\bigl\| X_{t}^{e_{j}} - X_{t}^{e_{k}}
\bigr\|_{1} > n^{-k}\bigr] < 2n^{2k+1-d},
\]
and so taking a union bound and applying the inequality proved just above,
\[
P \Bigl[\sup_{v,w \in\Delta_{n}} \bigl\| X_{t}^{v}
- X_{t}^{w} \bigr\|_{1} > n^{3-k}
\Bigr] < 2n^{2k+3-d}.
\]

This tells us that after $O(n \log n)$ steps, the $L^{1}$ distance
between any pair of points is extremely small with high probability.
The second\vspace*{1pt} step of the coupling is almost identical to the algorithm
given in\vspace*{1pt} Section \ref{sec4} of this note. Run $X_{t}^{(n^{-1},\ldots, n^{-1})}$
from time $T_{1}$ to time $T$, recording all choices of $i(t), j(t)$
and $\lambda(t)$ from representation (\ref{EqMoveRep}). Then form the
same partition process, and use it to attempt subset couplings of all
variables to this special chain. We will perform these couplings in
such a way that with high probability, all chains simultaneously have
successful subset couplings, rather than merely having a high
probability of a substantial fraction of the subset couplings
succeeding.

At each subset coupling stage, use the update variable $\lambda(t)$ for
the chain $X_{t}^{(n^{-1},\ldots, n^{-1})}$. For each other chain
$X_{t}^{v}$, there will be some probability $p(t,v)$ that $X_{t}^{v}$
performs a successful subset coupling with $X_{t}^{(n^{-1},\ldots,
n^{-1})}$. Let $p$ be a known lower bound on $\inf_{v \in\Delta_{n}}
p(t,v)$. This can be obtained from Lemma \ref{lem4.4} and inequality
(\ref{IneqCftp}). To determine the update value of $X_{t}^{v}$, choose
a single uniform random variable $U$. If $U < p$, let $X_{t}^{v}$ have
a successful subset coupling, in which case the change to $X_{t+1}^{v}$
depends only on $i(t),j(t)$ and $X_{t+1}^{(n^{-1},\ldots, n^{-1})}$,
not the particular value of $U$. Otherwise, update with $\lambda$ taken
from the $\frac{U - p}{1 - p}$'th quantile of the remainder
distribution. When $\frac{X_{t}^{v}[i] +
X_{t}^{v}[j]}{X_{t}^{(n^{-1},\ldots, n^{-1})}[i] +
X_{t}^{(n^{-1},\ldots, n^{-1})}[j]} \leq1$, this has density
\[
f(\lambda) = C \biggl( 1 - \frac{X_{t}^{v}[i] +
X_{t}^{v}[j]}{X_{t}^{(n^{-1},\ldots, n^{-1})}[i] +
X_{t}^{(n^{-1},\ldots, n^{-1})}[j]} \mathbf{1}_{g^{-1}(\lambda)
\in[0,1]} \biggr)
\]
for
\begin{eqnarray*}
g(\lambda) &=& \lambda\frac{X_{t}^{(n^{-1},\ldots, n^{-1})}[i] +
X_{t}^{(n^{-1},\ldots, n^{-1})}[j]}{X_{t}^{v}[i] + X_{t}^{v}[j]} \\
&&{}+
\frac{1}{X_{t}^{v}[i] + X_{t}^{v}[j]} \sum
_{s \in S \textbackslash\{ i
\} } \bigl(X_{t}^{(n^{-1},\ldots, n^{-1})}[s] -
X_{t}^{v}[s]\bigr)
\end{eqnarray*}
and $C$ a normalizing constant. An analogous formula holds when
\[
\frac{X_{t}^{v}[i] + X_{t}^{v}[j]}{X_{t}^{(n^{-1},\ldots, n^{-1})}[i]
+ X_{t}^{(n^{-1},\ldots, n^{-1})}[j]} > 1.
\]

Under this grand coupling, all subset couplings succeed together with
probability at least $p$. As long\vspace*{-1pt} as the $n$ points $X_{T_{1}}^{e_{j}}$
are close as measured in $L^{1}$ metric, and $X_{t}^{(n^{-1},\ldots,
n^{-1})}[i]$ remains far from 1 and 0, the proof of Lemma \ref{lem4.4}
tells us that all of the subset couplings succeed with high
probability. Finally, if a single subset coupling fails at time $t$,
then all chains should be coupled according to the proportional
coupling for time $s > t$.

It remains to determine what to do if one of the above subset couplings
fails. In order to obtain a perfect sample, it will be necessary to
look at a grand coupling for the epoch $-2T \leq t \leq-T$. Assume for
now that the grand coupling described above succeeds for the chain
started at $-2T$. It is necessary to determine the value at time 0 of
the chain started from $(\frac{1}{n},\ldots, \frac{1}{n})$ at time
$-2T$. Assume that at time~$-T$, this chain is at $v \in\Delta_{n}$.
Then our sample will be $X_{0}^{v}$. Fortunately, from the above
description, it is possible to calculate this value from $v$ and the
values of $i,j,\lambda$ and $U$ used during the first epoch. Thus, it
is sufficient to record those $O(T)$ pieces of information in each
failed epoch. A longer discussion of this algorithm, with pseudocode,
may be found in \cite{Smit12}.

It should be noted that, for other target distributions on the simplex,
such as those in Section \ref{sec6}, the above algorithm can also be used
without a rigorous bound on the mixing time and can be used to
rigorously check an estimated bound of time $T$. Simply run the
algorithm with epoch size $T$; the number of failed runs $k$ out of a
total of $N$ runs is distributed as a binomial random variable with
some unknown probability $q$, where $q$ is an upper bound on the total
variation distance to stationarity at time $T$.

\section*{Acknowledgments}

The author thanks David Aldous for mentioning the problem, and Olena
Blumberg, Persi Diaconis, Bob Hough, Daniel Jerison and John Jiang for
many helpful conversations. The author also thanks the reviewers for
friendly and useful comments.



\printaddresses

\end{document}